\documentclass[12pt,a4paper,leqno]{amsart}
\usepackage[english]{babel}
\usepackage[T1]{fontenc}
\usepackage{amssymb}
\usepackage{amsfonts}
\usepackage{amsmath}
\usepackage{mathptm}
\swapnumbers
\usepackage[english]{babel}
\usepackage[T1]{fontenc}
\usepackage{amssymb}
\usepackage{amsfonts}
\usepackage{amsmath}
\usepackage{mathptm}

\usepackage{eucal}
\begin{document}
\title{An inequality type condition for quasinearly subharmonic functions and applications}
\author{Juhani Riihentaus}
\vspace{3ex}
\noindent \address{{\quad}\newline
\noindent{University of Oulu} \newline
\noindent Department of Mathematical Sciences \newline \noindent P.O. Box 3000\newline
\noindent FI-90014 Oulun yliopisto \newline
\noindent Finland\newline
\noindent{and}\newline
\noindent{University of Eastern Finland}\newline
\noindent{Department of Physics and Mathematics}\newline \noindent{P.O. Box 111}\newline
\noindent{FI-80101 Joensuu}\newline
\noindent{Finland}}
\vspace{3ex}
\noindent\email{\noindent{riihentaus@member.ams.org and juhani.riihentaus@uef.fi}}
\begin{abstract} Generalizing older works of Domar and Armitage and Gardiner, we give an inequality for quasinearly subharmonic functions. As an application of this inequality, we improve   Domar's, Rippon's and our previous results concerning the existence of the largest subharmonic minorant of a given function. Moreover, and as an another application, we   give  a sufficient  condition for a separately quasinearly  subharmonic function to be quasinearly subharmonic. Our result  contains the previous results of Lelong, of Avanissian, of Arsove, of Armitage and Gardiner, and of ours.   
\end{abstract}
\date{September 19, 2015}
\maketitle
\section{Introduction}
\subsection{}  Let $D$ be a domain in ${\mathbb{R}^N}$, $N\geq 2$, and let $u:D\mapsto [-\infty ,+\infty )$ be subharmonic. We consider the function $u\equiv -\infty$ subharmonic.  Then $u$ is upper semicontinuous and satisfies the mean value inequality 
\begin{equation*}u(x)\leq \frac{1}{\nu _N\, r^N}\int\limits_{B^N(x,r)}u(y)\, dm_N(y)\end{equation*} 
for all balls $\overline{B^N(x,r)}\subset D$. It is an important fact that if $u$ is also nonnegative and $p>0$, then there exists a constant $C=C(N,p)$ such that
\begin{equation}u(x)^p\leq \frac{C}{\nu _N\, r^N}\int\limits_{B^N(x,r)}u(y)^p\, dm_N(y).\end{equation}

As a matter of fact, Fefferman and Stein \cite{FeSt72}, Lemma~2, p.~172, proved this inequality for absolute values of harmonic functions. See also \cite{Ga07}, Lemma~3.7, p.~116, \cite{Ku74}, Theorem~1, p.~529, and  \cite{AhRu93}, (1.5), p.~210 (also all these authors considered only absolute values of harmonic functions). However, the proof of Fefferman and Stein apply verbatim also in the more general case of nonnegative subharmonic functions. See \cite{Rii89}, Lemma, p.~69, and also \cite{Rii01}, \cite{Rii07}, \cite{Rii09}  and the references therein. A possibility for an essentially different proof was pointed out already in \cite{To86}, pp.~188-190. Later other different proofs were given in \cite{Pa94}, p.~18, and Theorem~1, p.~19 (see also \cite{Pa96}, Theorem~A, p.~15), \cite{Rii00}, Lemma~2.1, p.~233, and \cite{Rii01}, Theorem, p.~188.  Observe that the results in \cite{Pa94}, \cite{Rii00} and \cite{Rii01} hold, in addition to nonnegative subharmonic functions, also for more general quasinearly subharmonic functions. 

The inequality (1) has many applications. Among others, it has been applied to the weighted boundary behavior of subharmonic functions,  to the nonintegrability of subharmonic and superharmonic functions, and to the subharmonicity of separately subharmonic functions, see \cite{Rii89},  \cite{CeSa93}, \cite{Rii07}, \cite{Rii09}, \cite{DoRi13}, \cite{Rii13}, and the references therein.
\subsection{} In order to improve the above referred results on the subharmonicity of separately subharmonic functions, we give below in section~2 a rather general inequality type result which is related to the inequality (1), at least partly, and which applies more generally also to quasinearly subharmonic functions. This result has its origin in the previous considerations of Armitage and Gardiner \cite{ArGa93}, proof of Proposition~2, pp.~257-259, and in \cite{Rii08}, Lemma~3.2, p.~5, \cite{Rii11}, Lemma~2.2, p.~6. Observe, however and as already Armitage and Gardiner have pointed out, that this inequality is based on an old argument of Domar \cite{Do57}, proof of Proposition~2, pp.~257-259, and  proof of Theorem~1, pp.~258-259. In section~3 we will then apply the obtained inequality type result to domination conditions for families of quasinearly subharmonic functions, improving Domar's, Rippon's and our previous results, see \cite{Do57}, Theorem~1 and Theorem~2, pp.~430-431, \cite{Rip81}, Theorem~5, p.~128, and \cite{Rii11}, Theorem~2.1, pp.~4-5. In addition, in section~4  we apply this inequality to the quasinearly subharmonicity of separately quasinearly  subharmonic functions, slightly improving our previous results \cite{Rii08}, Theorem~4.1, pp.~8-9, \cite{Rii09}, Theorem~3.3.1, pp.~e2621-e2622. 

Though we indeed give improvements to our previous results, our presentation is, nevertheless and at least in some sense, of a survey type. Our notation is rather standard, see e.g. \cite{He71}, \cite{Rii07}, \cite{Rii08},  \cite{Rii09} and \cite{Rii11}. However and for the convenience of the reader, we begin by recalling  here the definitions of nearly subharmonic functions and quasinearly subharmonic functions. 
\subsection{Nearly subharmonic functions and quasinearly subharmonic functions}
We say that a function 
$u:\, D\rightarrow [-\infty ,+\infty )$ is \emph{nearly subharmonic}, if $u$ is Lebesgue measurable, $u^+\in {\mathcal{L}}^1_{\textrm{loc}}(D)$, 
and for all $\overline{B^N(x,r)}\subset D$,  
\begin{equation*}u(x)\leq \frac{1}{\nu _N\, r^N}\int\limits_{B^N(x,r)}u(y)\, dm_N(y).\end{equation*}
Observe that in the standard definition of nearly subharmonic functions one uses the slightly stronger  assumption that 
$u\in {\mathcal{L}}^1_{\textrm{loc}}(D)$, see e.g. \cite{He71}, p.~14. However, our above, slightly 
more general definition seems to be  more practical, see  e.g. \cite{Rii07},  Proposition~2.1~(iii) and Proposition~2.2~(vi), (vii), pp.~54-55, and \cite{Rii09},  Proposition~1.5.1~(iii) and Proposition~1.5.2~(vi), (vii), p.~e2615. The following lemma emphasizes this fact still more:

\vspace{2ex}

\noindent{\textbf{Lemma~1.1.}} (\cite{Rii07}, Lemma, p.~52) \textit{Let $u:\,D\rightarrow [-\infty ,+\infty )$ be Lebesgue measurable. 
Then  $u$ is nearly subharmonic (in the sense defined above) if and only if there exists a function $u^*$, subharmonic in $D$  such that $u^*\geq u$ and $u^*=u$ almost everywhere in $D$. Here $u^*$ is the upper semicontinuous regularization of $u$:
\begin{displaymath}u^*(x)=\limsup_{x'\rightarrow x}u(x').\end{displaymath}}

\vspace{2ex}

\noindent \textit{Proof}. 
The proof follows at once  from \cite{He71}, proof of Theorem~1, pp.~14-15  (and referring also to \cite{Rii07}, Proposition~2.1~(iii) and Proposition~2.2~(vii), pp.~54-55). \qed

\vspace{2ex}

We say that a Lebesgue measurable function $u:\,D \rightarrow 
[-\infty ,+\infty )$ is \emph{$K$-quasinearly subharmonic}, if  $u^+\in{\mathcal{L}}^{1}_{\textrm{loc}}(D)$ and if there is a 
constant $K=K(N,u,D)\geq 1$
such that for all  $\overline{B^N(x,r)}\subset D$,    
\begin{equation} u_M(x)\leq \frac{K}{\nu _N\,r^N}\int\limits_{B^N(x,r)}u_M(y)\, dm_N(y)\end{equation}
for all $M\geq 0$, where $u_M:=\max\{u,-M\}+M$. A function $u:\, D\rightarrow [-\infty ,+\infty )$ is \emph{quasinearly subharmonic}, if $u$ is 
$K$-quasinearly subharmonic for some $K\geq 1$.

\vspace{2ex}

For basic properties of quasinearly subharmonic functions,  see \cite{Rii07}, \cite{Rii09},  \cite{PaRi08}, and the references therein. We recall here only that this  function class 
includes, among \mbox{others,}  subharmonic functions, and, more generally,  quasisubharmonic and nearly subharmonic functions (see e.g. \cite{He71}, pp.~14, 26),    also functions satisfying certain natural  growth conditions, especially  
certain eigenfunctions, and  polyharmonic functions. Also, the class of Harnack functions is included, thus, among others, nonnegative harmonic functions 
as well as nonnegative solutions of some elliptic equations. In particular, the partial differential equations associated with quasiregular mappings 
belong to this family of elliptic equations. See e.g.~\cite{Vu82}.

Quasinearly subharmonic functions (perhaps with a different terminology, and sometimes in
 certain special cases, or just the corresponding generalized mean value inequalities (1) or (2)) have been considered in many papers, see e.g.  \cite{Rii89}, \cite{Pa94}, \cite{Rii07}, \cite{Rii08}, \cite{Rii09}, \cite{Rii11}, \cite{PaRi08}, \cite{Ko07}, \cite{DoRi10$_1$}, \cite{DoRi10$_2$}, \cite{KoMa12}, \cite{DoRi13},  and the  references therein. However and as a matter of fact,  already Domar \cite{Do57} considered (essentially) nonnegative quasinearly subharmonic functions.
\subsection{Two additional notational remarks}  The below presented proofs for our results, that is, the proofs of Theorem~2.1, Theorem~3.3 and Theorem~4.1, are quite much based on our previous arguments in \cite{Rii08}, proofs of Lemma~3.2, pp.~5-7,  and of Theorem~4.1, pp.~8-12, and \cite{Rii11}, proof of Theorem~2.1, pp.~4-8. Therefore, and in order to make  the possible comparison and checking etc. easier for the reader, we will  use, as previously, certain constants $s_0$, $s_1$, $s_2$, $s_3$, $s_4$ and $s_5$. 

Below in Examples 1, 2, 3, 4  and 5, we consider increasing functions   $\phi :[0,+\infty ]\rightarrow [0,+\infty ]$, say, which are of a certain form far away, that is, for big values of the argument. In such a case, we take the liberty to use the convention that the function is then automatically defined for small values of the argument in an approriate way.  As an example, if we write, for $p>0$, $q\in {\mathbb{R}}$, 
\[\phi (t)=\frac{t^p}{(\log t)^q},  \]
we mean the following function:
\begin{displaymath}
\varphi (t):=
\begin{cases}
\frac{t^p}{(\log t )^q}, &
 {\textrm{ when }}\, \,\,  t\geq t_1, \\
\frac{t}{t_1}\phi (t_1),
& {\textrm{ when }}\, 0\leq t<t_1,
\end{cases}
\end{displaymath}
where $t_1\geq 2$ is some suitable integer in  ${\mathbb{N}}$,  say.
\section{An inequality for quasinearly subharmonic functions} 
\subsection{} As pointed out already above, our previous result \cite{Rii08}, Lemma~3.2, p.~5, was a generalized version of Armitage's and Gardiner's argument \cite{ArGa93}, proof of  Proposition~2, pp.~257-258, and as such, it  was based on an old argument of Domar \cite{Do57}, Lemma~1, pp.~431-432 and 430. The following is another variant of this inequality type result:

\vspace{2ex}

\noindent{\textbf{Theorem~2.1.}}
\textit{Let $K\geq 1$.  Let $\varphi :[0,+\infty ]\rightarrow [0,+\infty ]$ and $\psi :[0,+\infty ]\rightarrow [0,+\infty ]$
 be increasing functions such that there are $s_0, \,s_1\in {\mathbb{N}}$,  $s_0<s_1$, such  that 
\begin{itemize}
\item[{(i)}] the inverse functions $\varphi ^{-1}$ and  $\psi^{-1}$ are defined on $[\min \{\,\varphi (s_1-s_0),\psi (s_1-s_0)\,\},+\infty ]$, 
\item[{(ii)}] $2K(\psi ^{-1}\circ \varphi )(s-s_0)\leq (\psi ^{-1}\circ \varphi )(s)$ for all $s\geq s_1$,
\item[{(iii)}] $\int\limits_{s_1}^{+\infty }\left[\frac{(\psi ^{-1}\circ \varphi )(s+2)}
{(\psi ^{-1}\circ \varphi )(s)}\frac{1}{\varphi (s-s_0)}\right]^{\frac{1}{N-1}}\, ds<+\infty .$
\end{itemize}
Let $u:\, D\rightarrow [0,+\infty )$ be a $K$-quasinearly subharmonic function. Let $\tilde {s}_1\in {\mathbb{N}}$, $\tilde {s}_1\geq s_1$, be arbitrary. 
Then for each $x\in D$ and $r>0$ such that 
$\overline{B^N(x,r)}\subset D$ either
\[u(x)\leq (\psi ^{-1}\circ \varphi )(\tilde {s}_1+1)\]
or
\[\Phi (u(x))\leq \frac{C}{r^N}\int\limits_{B^N(x,r)}\psi (u(y))\, dm_N(y)\]
where $C=C(N,K,s_0)$ and $\Phi :\, [0,+\infty )\rightarrow 
[0,+\infty ),$
}
\begin{displaymath}
 \Phi (t):=\begin{cases} \left( \int\limits_{(\varphi  ^{-1}\circ \psi )(t)-2}^{+\infty }
 \left[\frac{(\psi ^{-1}\circ \varphi )(s+2)}{(\psi ^{-1}\circ \varphi )(s)}\cdot \frac{1}{\varphi (s-s_0)}\right]^{\frac{1}{N-1}}\, ds \right)^{1-N},& 
{\textit{when }} t\geq s_3, \\
\frac{t}{s_3}\Phi (s_3), &{\textit{ when }\, } 0\leq t<s_3, 
\end{cases}
\end{displaymath} \textit{where}  $s_3:=\max\{\, s_1+3,s_2,(\psi ^{-1}\circ \varphi )(s_1+3)\,\}$ \textit{and}  $s_2:=\max\{\, s_1,(\psi ^{-1}\circ \varphi )(s_1+1)\,\}$.

\vspace{2ex}

\noindent{\textit{Proof}}.
The proof follows at once from \cite{Rii08}, proof of Lemma~3.2, pp.~5-7. As a matter of fact, it is sufficient to observe the following:
\begin{itemize}
\item Instead of functions $\varphi :[0,+\infty ) \rightarrow [0,+\infty )$ and $\psi :[0,+\infty ) \rightarrow [0,+\infty )$ one can equally well consider functions  $\varphi :[0,+\infty ] \rightarrow [0,+\infty ]$ and $\psi :[0,+\infty ] \rightarrow [0,+\infty ]$. See 
\cite{Rii11}, p.~5.
\item If
\[(\psi^{-1}\circ \varphi )(j_0)\leq t< (\psi^{-1}\circ \varphi )(j_0+1),\]
then 
\[\sum_{k=j_0}^{+\infty}\left[ \frac{(\psi^{-1}\circ \varphi )(k+1)}{(\psi^{-1}\circ \varphi )(k)}\cdot \frac{1}{\varphi (k-s_0)}\right]^{\frac{1}{N-1}}\leq 
\int\limits_{(\varphi^{-1}\circ \psi )(t)-2}^{+\infty}
\left[ \frac{(\psi^{-1}\circ \varphi )(s+2)}{(\psi^{-1}\circ \varphi )(s)}\cdot \frac{1}{\varphi (s-s_0)}\right]^{\frac{1}{N-1}}\, ds. 
\]
\end{itemize} \qed
\section{Domination conditions for families of quasinearly subharmonic functions}
\subsection{} We begin by recalling the results of Domar and Rippon. Let $F: D \rightarrow [0,+\infty ]$ be an upper semicontinuous function. Let ${\mathcal{F}}$ be a family of subharmonic functions $u: D\rightarrow [0,+\infty )$, which satisfy the condition
\[u(x)\leq F(x) \, {\textrm{ for all }}\,  x \in D.\]
Write 
\[w(x):=\sup_{u\in {\mathcal{F}}}\, u(x), \, \,\]
and let $w^*: D\rightarrow [0,+\infty ]$ be the upper semicontinuous regularization of $w$, that is,
\[w^*(x):=\limsup_{y\rightarrow x}\, w(y).\]
Improving the original results of Sjöberg \cite{Sj38} and Brelot \cite{Br45}, cf. also Green  \cite{Gr52}, Domar \cite{Do57}, Theorem~1 and Theorem~2, pp.~430-431, gave the following result:

\vspace{2ex}

\noindent{\textbf{Theorem~3.1.}}  \textit{If for some $\epsilon >0$,
\begin{equation*}\int\limits_D [\log^+F(x)]^{N-1+\epsilon}\, dm_N(x)<+\infty ,
\end{equation*}
then $w$ is locally bounded above in $D$, and thus $w^*$ is subharmonic in $D$.}

\vspace{2ex}

As Domar points out, his method of  proof applies also to   more general functions, that is,  to the above defined nonnegative quasinearly subharmonic functions. Much later Rippon \cite{Rip81},  Theorem~1,  p.~128, generalized   Domar's  result in the following form:

\vspace{2ex}

\noindent{\textbf{Theorem~3.2.}  \textit{Let $\varphi :[0,+\infty ]\rightarrow [0,+\infty ]$ be an increasing function such that 
\begin{equation*}\int\limits_1^{+\infty} \frac{dt}{\varphi (t)^{\frac{1}{N-1}}}<+\infty .
\end{equation*}
If 
\begin{equation*}\int\limits_D\varphi ( \log^+F(x))\, dm_N(x)<+\infty ,
\end{equation*}
then $w$ is locally bounded above in $D$, and thus $w^*$ is subharmonic in $D$.}

\vspace{2ex}

As pointed out by Domar, \cite{Do57}, pp.~436-440,  and by Rippon \cite{Rip81}, p.~129, the above results are for many particular cases sharp. 

In \cite{Rii11}, Theorem~2.1, pp.~4-5, we gave a general and at the same time flexible result which includes both Domar's and Rippon's results. Now we improve our result still further:

\vspace{2ex}

\noindent{\textbf{Theorem~3.3.}}  \textit{Let $K\geq 1$.  Let $\varphi :[0,+\infty ]\rightarrow [0,+\infty ]$ and $\psi :[0,+\infty ]\rightarrow [0,+\infty ]$
 be increasing functions for which there are $s_0, \,s_1\in {\mathbb{N}}$,  $s_0<s_1$, such  that 
\begin{itemize}
\item[{(i)}] the inverse functions $\varphi ^{-1}$ and  $\psi^{-1}$ are defined on $[\min \{\,\varphi (s_1-s_0),\psi (s_1-s_0)\,\},+\infty ]$, 
\item[{(ii)}] $2K(\psi ^{-1}\circ \varphi )(s-s_0)\leq (\psi ^{-1}\circ \varphi )(s)$ for all $s\geq s_1$,
\item[{(iii)}] the following integral is convergent:
\[ \int\limits_{s_1}^{+\infty }\left[\frac{(\psi ^{-1}\circ \varphi )(s+2)}{(\psi ^{-1}\circ \varphi )(s)}\cdot \frac{1}{\varphi (s-s_0)}\right]^{\frac{1}{N-1}}\,ds <+\infty .\]
\end{itemize}
Let ${\mathcal{F}}_K$ be a family of $K$-quasinearly subharmonic functions $u:\, D\rightarrow [-\infty ,+\infty )$ such that
\[u(x)\leq F_K(x) \,\,  for\, \,\,  all  \,\,\,  x\in D,\] 
where $F_K:\, D\rightarrow [0,+\infty ]$ is a Lebesgue measurable function. If for each compact set $E\subset D$,
\[\int\limits_E\psi (F_K(x))\, dm_N(x)<+\infty ,\]
then the family ${\mathcal{F}}_K$ is locally (uniformly) bounded in $D$. Moreover, the function $w^*:\, D\rightarrow [0,+\infty )$ is  \mbox{$K$-quasinearly} subharmonic. Here
\[w^*(x):=\limsup_{y\rightarrow x}w(y),\]
where \[w(x):=\sup_{u\in {\mathcal{F}}_K}u^+(x).\]
}

\vspace{2ex}

\noindent{\textit{Proof}}. 
Let $E$ be an arbitrary compact subset of $D$. Write $\rho_0:=$dist$(E,\partial D)$. Clearly $\rho_0>0$. Write 
\begin{displaymath}E_1:=\bigcup_{x\in E}\overline{B^n(x,\frac{\rho_0}{2})}.\end{displaymath}
Then $E_1$ is compact, and $E\subset E_1\subset D$.  Take $u\in {\mathcal{F}}_K^+$ arbitrarily,
where
\[{\mathcal{F}}_K^+:=\{\,u^+\, :\, u\in {\mathcal{F}}_K\, \}.\]
Let $\tilde{s}_1=s_1+2$, say. Take $x\in E$ arbitrarily and suppose that $u(x)> \tilde{s}_3$, where \newline \mbox{$\tilde{s}_3:=\max \{\, s_1+5,(\psi^{-1}\circ \varphi )(s_1+5)\,\}$}, say. By Theorem~2.1 we have,
\begin{equation*}\Phi (u(x)) 
\leq \frac{C}{\left(\frac{\rho_0}{2}\right)^N}\int\limits_{B^N(x,\frac{\rho_0}{2})}\psi (u(y))\, dm_N(y)
\leq \frac{C}{\left(\frac{\rho_0}{2}\right)^N}\int\limits_{E_1}\psi (F_K(y))\, dm_N(y)<+\infty ,
\end{equation*}
where
\[\Phi (t):=\left( \int\limits_{(\varphi^{-1}\circ \psi )(t)-2}^{+\infty}\left[\frac{(\psi ^{-1}\circ \varphi )(s+2)}{(\psi ^{-1}\circ \varphi )(s)}\cdot \frac{1}{\varphi (s-s_0)}\right]^{\frac{1}{N-1}}\, ds \right)^{1-N}.\]
Now we know that 
\[\int\limits_{s_1}^{+\infty}\left[\frac{(\psi ^{-1}\circ \varphi )(s+2)}{(\psi ^{-1}\circ \varphi )(s)}\cdot \frac{1}{\varphi (s-s_0)}\right]^{\frac{1}{N-1}}\, ds<+\infty .\]
Therefore, the set of function values 
\begin{equation*}
u(x), \quad x\in E, \quad u\in {\mathcal{F}}_K^+,
\end{equation*}
is  bounded above.

The rest of the proof goes then as in \cite{Rii11}, proof of Theorem~2.1, pp.~7-8.
\qed

\vspace{2ex}

\noindent{\textit{Remark~3.4.}} In Theorem~3.3 one can, instead of the assumption (iii), use also the following:
\begin{itemize}
\item[{(iii')}] \emph{the following series is convergent:}
\[ \sum\limits_{j=s_1+1}^{+\infty }\left[\frac{(\psi ^{-1}\circ \varphi )(j+1)}{(\psi ^{-1}\circ \varphi )(j)}\cdot \frac{1}{\varphi (j-s_0)}\right]^{\frac{1}{N-1}}<+\infty .\]
\end{itemize}

Instead of the above used function $\Phi$ one uses then the function $\Phi_1: [s_2,+\infty )\rightarrow [0,+\infty )$,
\begin{displaymath}
 \Phi_1 (t):=\left( \sum\limits_{k=j_0}^{+\infty }
 \left[\frac{(\psi ^{-1}\circ \varphi )(k+1)}{(\psi ^{-1}\circ \varphi )(k)}\cdot \frac{1}{\varphi (k-s_0)}\right]^{\frac{1}{N-1}} \right)^{1-N},
 \end{displaymath}
where $j_0\in \{\, s_1,s_1+2, \dots \, \}$ is such that 
\[(\psi^{-1}\circ \varphi )(j_0)\leq t< (\psi^{-1}\circ \varphi )(j_0+1).\]
The function $\Phi_1$ may of course be extended to the whole interval $[0,+\infty )$, for example as follows:
\begin{displaymath}
 \Phi_1 (t):=\begin{cases} \Phi_1(t),& 
{\textrm{ when }} t\geq s_3, \\
\frac{t}{s_3}\Phi_1 (s_3), &{\textrm{ when }\, } 0\leq t<s_3. 
\end{cases}
\end{displaymath}

Before giving examples, we write Theorem~3.3 in the following, perhaps more concrete form,  see also \cite{Rii11}, Remark~2.5, p.~8.

\vspace{2ex}

\noindent{\textbf{Corollary~3.5.}}  \textit{Let $K\geq 1$.  Let $\varphi :[0,+\infty ]\rightarrow [0,+\infty ]$ and $\phi :[0,+\infty ]\rightarrow [0,+\infty ]$
 be strictly increasing surjections for which there are $s_0, \,s_1\in {\mathbb{N}}$,  $s_0<s_1$, such  that 
\begin{itemize}
\item[{(i)}] $2K\phi^{-1}(e^{s-s_0})\leq \phi ^{-1}(e^s)$ for all $s\geq s_1$,
\item[{(ii)}] the following integral is convergent:
\[ \int\limits_{s_1}^{+\infty }\left[\frac{\phi^{-1} (e^{s+2})}{\phi^{-1} (e^s)}\frac{1}{\varphi (s-s_0)}\right]^{\frac{1}{N-1}}\,ds<+\infty .\]
\end{itemize}
Let ${\mathcal{F}}_K$ be a family of $K$-quasinearly subharmonic functions $u:\, D\rightarrow [-\infty ,+\infty )$ such that
\[u(x)\leq F_K(x) \,\,\, for\,\, all \,\,\, x\in D,\]
where $F_K:\, D\rightarrow [0,+\infty ]$ is a Lebesgue measurable function. If for each compact set $E\subset D$,
\[\int\limits_E\varphi (\log^+\phi (F_K(x)))\, dm_N(x)<+\infty ,\]
then the family ${\mathcal{F}}_K$ is locally (uniformly) bounded in $D$. Moreover, the function $w^*:\, D\rightarrow [0,+\infty )$ is  \mbox{$K$-quasinearly} subharmonic. Here
\[w^*(x):=\limsup_{y\rightarrow x}w(y),\]
where \[w(x):=\sup_{u\in {\mathcal{F}}_K}u^+(x).\]
}

\vspace{2ex}

\noindent{\textit{Proof.}} For the proof, just  choose $\psi (t)=\varphi (\log^+\phi (t))$. Since only big values count, we may simply use the formula  $\psi (t)=\varphi (\log \phi (t))$. One sees easily that for some $\tilde{s}_2\geq s_1$,
\begin{displaymath}
(\psi^{-1}\circ \varphi )(s)=\phi^{-1}(e^s) \,\, {\textrm{ for all }} \,\, s\geq \tilde{s}_2 .
\end{displaymath}
It is then  easy to see that the assumptions of Theorem~3.3 are satisfied.
\qed

\vspace{2ex}

\noindent{\textit{Example~1.}}
Let $\varphi :\, [0, +\infty ]\rightarrow [0,+\infty ]$ be a strictly increasing surjection such that (for some $s_0\in {\mathbb{N}}$),
\begin{equation*}\int\limits_{s_1}^{+\infty} \frac{ds}{\varphi (s-s_0)^{\frac{1}{N-1}}}<+\infty .
\end{equation*}
Choosing then various functions $\phi$ which, together with $\varphi$, satisfy the conditions (i) and (ii) of Corollary~3.5, one gets more concrete results. If  $\phi$ and $\phi^{-1}$ satisfy (at least far away) the $\Delta_2$-condition, then the conditions (i) and  (ii) are surely satisfied (see also  \cite{Rii11}, Remark~2.5, p.~8). Typical choices for $\phi$ might be, say, the following: 
\[\phi (t):=\frac{t^p}{(\log t)^q}, \, \, p>0, \, q\in {\mathbb{R}}.\]
The choice $p=1$, $q=0$ gives then the results of Domar and Rippon, Theorem~3.1 and Theorem~3.2 above. Choosing $0<p<1$ and $q\geq 0$, one gets (at least formal) improvements.

\vspace{2ex}

\noindent{\textit{Example~2.}} Let $\phi :\, [0, +\infty ]\rightarrow [0,+\infty ]$ be a strictly increasing surjection for which there is $t_1>1$ such that 
\begin{equation*}s=\phi (t)=\log t \,\,\,\, {\textrm{}}{\textrm{for }}\,\, t\geq t_1 {\textrm{}} 
\end{equation*}
and
\begin{equation*}t=\phi^{-1} (s)=e^s \,\,\,\, {\textrm{}}{\textrm{for }}\,\, s\geq \log t_1{\textrm{}}. 
\end{equation*}
One sees easily that the condition (i) of Corollary~3.5 holds. In this case $\phi^{-1}$ does not satisfy a $\Delta_2$-condition. Therefore, in order to apply Corollary~3.5 to a family ${\mathcal{F}}_K$ of $K$-quasinearly subharmonic functions, we must choose an approriate strictly increasing surjection $\varphi :\, [0,+\infty ]\rightarrow [0,+\infty ]$ such that for some $\tilde{s}_1\in {\mathbb{N}}$,
\begin{equation*}
\int\limits_{\tilde{s}_1}^{+\infty} \left[\frac{\phi^{-1}(e^{s+2})}{\phi^{-1}(e^{s})}\cdot \frac{1}{\varphi (s-s_0)}\right]^{\frac{1}{N-1}}ds=\int\limits_{\tilde{s}_1}^{+\infty} \frac{e^{\frac{e^2-1}{N-1}e^s}}{\varphi (s-s_0)^{\frac{1}{N-1}}}ds <+\infty .
\end{equation*}
Therefore, we have two restrictions for $\varphi$. As a first condition the above quite strong restriction and, as a second one,  the following, at least seemingly mild condition: For each compact set $E\subset D$,
\begin{equation*}
\int\limits_{E} \varphi (\log^+ (\log F_K(x)))\, dm_N(x) <+\infty .
\end{equation*}
\section{On the subharmonicity of separately subharmonic functions and generalizations} 
\subsection{} Wiegerinck \cite{Wi88}, Theorem, p.~770, see also \cite{WiZe91}, Theorem~1, p.~246, has shown that  separately subharmonic functions need not be subharmonic. On the other hand, Armitage and Gardiner \cite{ArGa93}, Theorem~1, p.~256, showed that a separately subharmonic function $u$ in a domain $\Omega$ of ${\mathbb{R}}^{m+n}$, $m\geq n\geq 2$, is subharmonic, 
provided $\phi \circ \log^+u^+\in {\mathcal{L}}^1_{{\textrm{loc}}}(\Omega )$, where $\phi :\, [0,+\infty )\rightarrow [0,+\infty )$ is an increasing function such that 
\begin{displaymath}
\int\limits_{1}^{+\infty }\frac{s^{\frac{n-1}{m-1}}}{\phi (s)^{\frac{1}{m-1}}}ds<+\infty .
\end{displaymath}
Armitage's and Gardiner's result included all the  previous existing results, that is, the results  of  Lelong \cite{Le45},  Théorème~1~bis, p.~315, of Avanissian \cite{Av61}, Théorème~9, p.~140, see also \cite{Le69}, Proposition~3, p.~24,  and \cite{He71}, Theorem, p.~31, of  Arsove \cite{Ar66}, Theorem~1, p.~622,  and ours \cite{Rii89}, Theorem~1, p.~69. Though Armitage's and Gardiner's result was close to being sharp, see \cite{ArGa93}, p.~255, it was, however, possible to improve their result slightly further. This was done in \cite{Rii08}, Theorem~4.1, pp.~8-9, with the aid of quasinearly subharmonic functions. See also \cite{Rii09}, Theorem~3.3.1 and Corollary~3.3.3, pp.~e2621-e2622. Now we improve our result still further:

\vspace{2ex}

\noindent{\textbf{Theorem~4.1.}} \textit{Let $K\geq 1$. Let $\Omega $ be a domain in ${\mathbb{R}}^{m+n}$, \mbox{$m\geq n\geq 2$.}
Let $u:\, \Omega \rightarrow [-\infty ,+\infty )$ be  a Lebesgue measurable function. Suppose that the following conditions are satisfied:}
\begin{itemize}
\item[(a)] \textit{For each $y\in {\mathbb{R}}^n$ the function}
\[\Omega (y)\ni x\mapsto u(x,y)\in [-\infty ,+\infty )\]
\textit{is $K$-quasinearly subharmonic.}
 \item[(b)] \textit{For each $x\in {\mathbb{R}}^m$ the function}
\[\Omega (x)\ni y\mapsto u(x,y)\in [-\infty ,+\infty )\]
\textit{is $K$-quasinearly subharmonic.}
\item[(c)] \textit{There are increasing functions  $\varphi :[0,+\infty )\rightarrow [0,+\infty )$ and $\psi :[0,+\infty )\rightarrow [0,+\infty )$
 and  $s_0, \,s_1\in {\mathbb{N}}$,  $s_0<s_1$, such  that}
\begin{itemize}
\item[{(c1)}] \textit{the inverse functions $\varphi^{-1}$ and $\psi ^{-1}$ are defined on $[\min \{\, \varphi (s_1-s_0),\psi (s_1-s_0)\,\},+\infty )$,}
 \item[{(c2)}] \textit{$2K(\psi ^{-1}\circ \varphi )(s-s_0)\leq (\psi ^{-1}\circ \varphi )(s)$ for all $s\geq s_1$,}
\item[{(c3)}] \textit{the following integral is convergent:}
\[\int\limits_{s_5}^{+\infty }\left[\frac{(\psi ^{-1}\circ \varphi )(s+2)}{(\psi ^{-1}\circ \varphi )(s)}\cdot \frac{1}{\varphi (s-s_0)}\right]^{\frac{1}{m-1}}\cdot\left( \int\limits_{s_5+s_0+2}^{s+s_0+2}\left[ \frac{(\psi ^{-1}\circ \varphi )(t+2)}{(\psi ^{-1}\circ \varphi )(t)} \right]^{\frac{1}{n-1}}\, dt\right)^{\frac{n-1}{m-1}}\,ds<+\infty , \]
\textit{where} $s_5:=\max \{\, s_0+s_1+3, s_0+(\psi ^{-1}\circ \varphi )(s_1+3), s_0+(\varphi ^{-1}\circ \psi  )(s_1+3)\,\}$,
\item[{(c4)}] $\psi \circ u^+\in {\mathcal{L}}_{\textrm{loc}}^{1}(\Omega )$.
\end{itemize}
\end{itemize}
\textit{Then} $u$
\textit{is quasinearly subharmonic in $\Omega $}.

\vspace{2ex}

\noindent{\textit{Proof}}. Recall that $s_2=\max \{\,s_1, (\psi ^{-1}\circ \varphi )(s_1+1)\,\}$ and  $s_3=\max \{\,s_1+3,s_2, (\psi ^{-1}\circ \varphi )(s_1+3)\,\}$. Write $s_4:=\max \{\, s_0+s_3, (\varphi ^{-1}\circ \psi  )(s_1+3)\,\}$, say. Clearly, $s_0<s_1< s_3<s_4< s_5$. (We may of course suppose that   $s_3$, $s_4$ and $s_5$ are integers.) One may replace $u$ by $\max \{\, u^+,M\,\}$,
where $M=\max \{\, s_5+3,(\psi ^{-1}\circ \varphi )(s_4+3), (\varphi ^{-1}\circ \psi )(s_4+3)\,\}$, say. We  continue to denote $u_M$ by $u$.

\vspace{0.5ex}

\noindent{\textbf{Step 1}} \textit{Use of Theorem~2.1}

\vspace{0.5ex}

Take $(x_0,y_0)\in \Omega $ and $r>0$ arbitrarily such that $\overline{B^m(x_0,2r)\times B^n(y_0,2r)}\subset \Omega $. Take $(\xi , \eta )\in B^m(x_0,r)\times B^n(y_0,r)$ arbitrarily.  We know that $u(\cdot ,y)$ is \mbox{$K$-quasinearly} subharmonic for each $y\in B^n(y_0,2r)$. In order to apply Theorem~2.1, it is clearly sufficient to show that  
\[\int\limits_{s_5+s_0+2}^{+\infty}\left[\frac{(\psi^{-1}\circ \varphi )(s+2)}{(\psi^{-1}\circ \varphi )(s)}\cdot \frac{1}{\varphi (s-s_0)}\right]^{\frac{1}{m-1}}\, ds<+\infty .\]
But this follows at once from the assumption (c3), since for all $s\geq s_5+s_0+2$,
\[\left( \int\limits_{s_5+s_0+2}^{s+s_0+2}\left[\frac{(\psi^{-1}\circ \varphi )(t+2)}{(\psi^{-1}\circ \varphi )(t)}\right]^{\frac{1}{n-1}}\, dt\right)^{\frac{n-1}{m-1}}\geq \left(\int\limits_{s_5+s_0+2}^{s+s_0+2}1\, dt\right)^{\frac{n-1}{m-1}}=(s-s_5)^{\frac{n-1}{m-1}}\geq (s_0+2)^{\frac{n-1}{m-1}} .\]
From Theorem~2.1 it then follows that for all $y\in B^n(y_0,2r)$
\begin{equation}
\tilde{\Phi}(u(\xi ,y))\leq \frac{C}{r^m}\int\limits_{B^m(\xi ,r)}\psi (u(x,y))dm_m(x),
\end{equation}
where
\begin{displaymath}
\tilde{\Phi}(t):=\begin{cases}\left(\int\limits^{+\infty }_{(\varphi ^{-1}\circ \psi )(t) -2}\left[\frac{(\psi^{-1}\circ \varphi )(s+2)}{(\psi^{-1}\circ \varphi )(s)}\frac{1}{\varphi (s-s_0)}\right]^{\frac{1}{m-1}}ds\right)^{1-m}, & {\textrm{ when }}\, t\geq s_3,\\
\frac{t}{s_3}\tilde{\Phi}(s_3),
& {\textrm{when }}\, 0\leq t<
s_3.
\end{cases}\end{displaymath}

\vspace{0.5ex}

\noindent{\textbf{Step 2}} \textit{Take mean values on both sides of} (3)

\vspace{0.5ex}

Taking  (generalized) mean values with respect to the variable $y$ over the ball $B^n(\eta ,r)$ on both sides of (3), we get:
\begin{equation*}\begin{split}
\frac{C}{r^n}\int\limits_{B^n(\eta ,r)}\tilde{\Phi}(u(\xi ,y))dm_n(y)
&\leq \frac{C}{r^n}\int\limits_{B^n(\eta ,r)}[\frac{C}{r^m}\int\limits_{B^m(\xi ,r)}
\psi (u(x,y))dm_m(x)]\, dm_n(y)\\
&\leq \frac{C}{r^{m+n}}\int\limits_{B^m(\xi ,r)\times B^n(\eta ,r)}
\psi (u(x,y))dm_{m+n}(x,y)\\
&\leq \frac{C}{r^{m+n}}\int\limits_{B^m(x_0,2r)\times B^n(y_0,2r)}
\psi (u(x,y))dm_{m+n}(x,y).
\end{split}
\end{equation*}
Here one must of course check that both $\psi \circ u(\cdot ,\cdot )$ and $\tilde{\Phi}(u(\xi ,\cdot ))$ are  Lebesgue measurable!

\vspace{0.5ex}

\noindent{\textbf{Step 3}} \textit{In order to apply Theorem~2.1 once more, define new functions $\varphi_1$ and $\psi_1$}

\vspace{0.5ex}

Write $\psi_1:\, [0,+\infty )\rightarrow [0,+\infty )$,
\begin{displaymath}
\psi_1(t):=\tilde{\Phi}(t)=\begin{cases}\left(\int\limits^{+\infty }_{(\varphi ^{-1}\circ \psi )(t) -2}\left[\frac{(\psi^{-1}\circ \varphi )(s+2)}{(\psi^{-1}\circ \varphi )(s)}\frac{1}{\varphi (s-s_0)}\right]^{\frac{1}{m-1}}ds\right)^{1-m}, & {\textrm{ when }}\, t\geq s_3,\\
\frac{t}{s_3}\tilde{\Phi}(s_3),
& {\textrm{when }}\, 0\leq t<
s_3.
\end{cases}\end{displaymath}
It is easy to see that $\psi_1$ is defined, strictly increasing and continuous.
Write then $\varphi_1 :\, [0,+\infty )\rightarrow [0,+\infty )$,
\begin{displaymath}
\varphi _1(t):=\begin{cases}  \psi _1((\psi ^{-1}\circ \varphi )(t))=\tilde{\Phi}(\psi ^{-1}(\varphi (t))), & {\textrm{when }}\, t\geq s_3,\\
\frac{t}{s_3}\psi_1 ((\psi ^{-1}\circ \varphi )(s_3))=\frac{t}{s_3}\tilde{\Phi}(\psi ^{-1}(\varphi (s_3))),
& {\textrm{when }}\, 0\leq t<s_3. \\
\end{cases}\end{displaymath}
Also $\varphi_1$ is defined, strictly increasing and continuous. This follows from the facts that $\psi_1$ is defined, strictly increasing and continuous (similarly as the functions $\varphi \vert [s_1-s_0,+\infty)$ and  $\psi \vert [s_1-s_0,+\infty)$).
Observe here that for $t\geq s_4$, say,
\begin{equation*}\begin{split}
\varphi _1(t)&=
\left(\int\limits_{(\varphi ^{-1}\circ \psi )((\psi ^{-1}\circ \varphi )(t))-2}^{+\infty }
\left[\frac{(\psi^{-1}\circ \varphi )(s+2)}{(\psi^{-1}\circ \varphi )(s)}\frac{1}{\varphi (s-s_0)}\right]^{\frac{1}{m-1}}\, ds\right)^{1-m}\\
&=\left(\int\limits_{t-2}^{+\infty }\left[\frac{(\psi^{-1}\circ \varphi )(s+2)}{(\psi^{-1}\circ \varphi )(s)}\frac{1}{\varphi (s-s_0)}\right]^{\frac{1}{m-1}}\, ds \right)^{1-m}.\end{split}\end{equation*}
One sees easily that $(\psi_1^{-1}\circ \varphi_1)(t)=(\psi^{-1}\circ \varphi)(t)$ for all $t\geq s_3$, thus $2K(\psi_1^{-1}\circ \varphi_1)(s-s_0)\leq (\psi_1^{-1}\circ \varphi_1)(s)$  for all $s\geq \overline{s}_1\geq s_4$.

To show that
\[
\int\limits_{s_5+s_0+2}^{+\infty }\left[\frac{(\psi_1^{-1}\circ \varphi_1 )(s+2)}{(\psi_1^{-1}\circ \varphi_1 )(s)}\frac{1}{\varphi_1 (s-s_0)}\right]^{\frac{1}{n-1}}\, ds
 <+\infty ,\]
we proceed as follows.   

Write 
$F:\, [s_5,+\infty )\times [s_5+s_0+2,+\infty )\rightarrow [0,+\infty )$,
\begin{displaymath}
F(s,t):=\begin{cases}
\left[\frac{(\psi_1^{-1}\circ \varphi_1)(t+2)}
{(\psi_1^{-1}\circ \varphi_1)(t)}\frac{(\psi^{-1}\circ \varphi)(s+2)}
{(\psi^{-1}\circ \varphi)(s)}\frac{1}{\varphi(s-s_0)}\right]^{\frac{1}{m-1}}, & {\textrm{when }}\, s_5+s_0+2\leq t-s_0-2\leq s,\\
0,
& {\textrm{when }}\, s_5\leq s<t-s_0-2. 
\end{cases}\end{displaymath}
Suppose that $m>n\geq 2$. Then just calculate, use Minkowski's inequality and assumption (c3): 
\begin{align*}
&\left(\int\limits_{s_5+s_0+2}^{+\infty}\left[\frac{(\psi_1^{-1}\circ \varphi_1)(t+2)}{(\psi_1^{-1}\circ \varphi_1)(t)}\frac{1}{\varphi_1(t-s_0)}\right]^{\frac{1}{n-1}}dt\right)^{\frac{n-1}{m-1}}=
\\
&=\left(\int\limits_{s_5+s_0+2}^{+\infty}\left[\frac{(\psi_1^{-1}\circ \varphi_1)(t+2)}
{(\psi_1^{-1}\circ \varphi_1)(t)}\right]^{\frac{1}{n-1}}
\left( \int\limits_{t-s_0-2}^{+\infty}\left[ \frac{(\psi^{-1}\circ \varphi)(s+2)}
{(\psi^{-1}\circ \varphi)(s)}\frac{1}{\varphi(s-s_0)}\right]^{\frac{1}{m-1}}ds 
\right)^{-\frac{1-m}{n-1}}dt\right)^{\frac{n-1}{m-1}}\\
&=\left(\int\limits_{s_5+s_0+2}^{+\infty}\left[\frac{(\psi_1^{-1}\circ \varphi_1)(t+2)}
{(\psi_1^{-1}\circ \varphi_1)(t)}\right]^{\frac{1}{n-1}}
\left( \int\limits_{t-s_0-2}^{+\infty}\left[ \frac{(\psi^{-1}\circ \varphi)(s+2)}
{(\psi^{-1}\circ \varphi)(s)}\frac{1}{\varphi(s-s_0)}\right]^{\frac{1}{m-1}}ds 
\right)^{\frac{m-1}{n-1}}dt\right)^{\frac{n-1}{m-1}}\\
&=\left(\int\limits_{s_5+s_0+2}^{+\infty}
\left(\int\limits_{t-s_0-2}^{+\infty}\left(\left[\frac{(\psi_1^{-1}\circ \varphi_1)(t+2)}
{(\psi_1^{-1}\circ \varphi_1)(t)}\cdot \frac{(\psi^{-1}\circ \varphi)(s+2)}
{(\psi^{-1}\circ \varphi)(s)}\frac{1}{\varphi(s-s_0)}\right]^{\frac{1}{m-1}}\right)ds 
\right)^{\frac{m-1}{n-1}}dt\right)^{\frac{n-1}{m-1}}\\
&=\left( \int\limits_{s_5+s_0+2}^{+\infty}\left[\int\limits_{s_5}^{+\infty}F(s,t)\, ds\right]^{\frac{m-1}{n-1}}dt\right)^{\frac{n-1}{m-1}}
\leq \left( \int\limits_{s_5}^{+\infty}\left[\int\limits_{s_5+s_0+2}^{+\infty}F(s,t)^{\frac{m-1}{n-1}}\, dt\right]^{\frac{n-1}{m-1}}ds\right)=\\
&=\int\limits_{s_5}^{+\infty}
\left(\int\limits_{s_5+s_0+2}^{s+s_0+2}\left(\left[\frac{(\psi_1^{-1}\circ \varphi_1)(t+2)}
{(\psi_1^{-1}\circ \varphi_1)(t)}\right]^{\frac{1}{m-1}}\left[ \frac{(\psi^{-1}\circ 
\varphi)(s+2)}
{(\psi^{-1}\circ \varphi)(s)}\frac{1}{\varphi(s-s_0)}\right]^{\frac{1}{m-1}}
\right)^{\frac{m-1}{n-1}}dt 
\right)^{\frac{n-1}{m-1}}ds\\
&=\int\limits_{s_5}^{+\infty}\left[ \frac{(\psi^{-1}\circ 
\varphi)(s+2)}
{(\psi^{-1}\circ \varphi)(s)}\frac{1}{\varphi(s-s_0)}\right]^{\frac{1}{m-1}}
\left(\int\limits_{s_5+s_0+2}^{s+s_0+2}\left[ \frac{(\psi_1^{-1}\circ \varphi_1)(t+2)}
{(\psi_1^{-1}\circ \varphi_1)(t)}\right]^{\frac{1}{n-1}}dt\right)^{\frac{n-1}{m-1}}ds \\ 
&=\int\limits_{s_5}^{+\infty}\left[ \frac{(\psi^{-1}\circ 
\varphi)(s+2)}
{(\psi^{-1}\circ \varphi)(s)}\frac{1}{\varphi(s-s_0)}\right]^{\frac{1}{m-1}}
\left(\int\limits_{s_5+s_0+2}^{s+s_0+2}\left[ \frac{(\psi^{-1}\circ \varphi)(t+2)}
{(\psi^{-1}\circ \varphi)(t)}\right]^{\frac{1}{n-1}}dt\right)^{\frac{n-1}{m-1}}ds <+\infty .\end{align*}
The case $m=n$ is considered similarly,  just replacing  Minkowski's inequality with Fubini's theorem.

\vspace{0.5ex}

\noindent{\textbf{Step 4}} \textit{Apply Theorem~2.1 to conclude that $u(\cdot ,\cdot )$ is bounded in $B^m(x_0,r)\times B^n(y_0,r)$}

\vspace{0.5ex}

With the aid of Theorem~2.1 we get 
\begin{equation*}
\Psi(u(\xi ,\eta ))
\leq \frac{C}{r^n}\int\limits_{B^n(\eta  ,r)}\tilde{\Phi} (u(\xi ,y))dm_n(y)
\leq \frac{C}{r^{m+n}}\int\limits_{B^m(x_0,2r)\times B^n(y_0,2r)}\psi (u(x,y))dm_{m+n}(x,y),
\end{equation*}
where now 
\begin{displaymath}
\Psi (t):=\begin{cases}\left( \int\limits_{(\varphi_1^{-1}\circ \psi_1 )(t)-2}^{+\infty}\left[ \frac{(\psi^{-1}_1\circ  \varphi_1)(s+2)} {(\psi^{-1}_1\circ  \varphi_1)(s)}\frac{1}{\varphi_1(s-s_0)}\right]^{\frac{1}{n-1}}ds \right)^{1-n},
& {\textrm{when }}\, t\geq s_3, \\
 \frac{t}{s_3}\Psi (s_3), & {\textrm{when }}\, 0\leq t< s_3,
\end{cases}
\end{displaymath}
or equivalently
\begin{displaymath}
\Psi (t):=\begin{cases}\left( \int\limits_{(\varphi^{-1}\circ \psi )(t)-2}^{+\infty}\left[ \frac{(\psi^{-1}_1\circ  \varphi_1)(s+2)} {(\psi^{-1}_1\circ  \varphi_1)(s)}\frac{1}{\varphi_1(s-s_0)}\right]^{\frac{1}{n-1}}ds \right)^{1-n},
& {\textrm{when }}\, t\geq s_3, \\
 \frac{t}{s_3}\Psi (s_3), & {\textrm{when }}\,  0\leq t< s_3.
\end{cases}
\end{displaymath}
Observe that we know that
\begin{displaymath}
\begin{split}
\int\limits_{s_5+s_0+2}^{+\infty}& \left[ \frac{(\psi^{-1}_1\circ  \varphi_1)(s+2)} {(\psi^{-1}_1\circ  \varphi_1)(s)}\frac{1}{\varphi_1(s-s_0)}\right]^{\frac{1}{n-1}}ds\leq \\ 
&
\leq 
\int\limits_{s_5}^{+\infty}\left[\frac{(\psi^{-1}\circ \varphi)(s+2)}
{(\psi^{-1}\circ \varphi)(s)}\frac{1}{\varphi(s-s_0)}\right]^{\frac{1}{m-1}}
\left( \int\limits_{s_5+s_0+2}^{s+s_0+2}\left[ \frac{(\psi^{-1}\circ \varphi)(t+2)}
{(\psi^{-1}\circ \varphi)(t)}\right]^{\frac{1}{n-1}}dt 
\right)^{\frac{n-1}{m-1}}ds,
 \end{split}
\end{displaymath}
and that by assumption (c3),
\begin{displaymath}
\int\limits_{s_5}^{+\infty}\left[\frac{(\psi^{-1}\circ \varphi)(s+2)}
{(\psi^{-1}\circ \varphi)(s)}\frac{1}{\varphi(s-s_0)}\right]^{\frac{1}{m-1}}
\left( \int\limits_{s_5+s_0+2}^{s+s_0+2}\left[ \frac{(\psi^{-1}\circ \varphi)(t+2)}
{(\psi^{-1}\circ \varphi)(t)}\right]^{\frac{1}{n-1}}dt 
\right)^{\frac{n-1}{m-1}}ds<+\infty .
\end{displaymath}
Hence the set of function values
\[(\varphi_1^{-1}\circ \psi_1 )(u(\xi ,\eta ))-2=(\varphi^{-1}\circ \psi )(u(\xi ,\eta ))-2, \,\, (\xi ,\eta )\in B^m(x_0,r)\times B^n(y_0,r),\]
must be bounded. Thus the function $u(\cdot ,\cdot )$ is bounded above in  $B^m(x_0,r)\times B^n(y_0,r)$. By \cite{Rii07}, Proposition~3.1, p.~57, (or by \cite{Rii09}, Proposition~3.2.1, p.~e2620)  we see that $u(\cdot ,\cdot )$ is quasinearly subharmonic. 
\qed

\vspace{2ex}

\noindent{\textbf{Corollary~4.2.}}
 \textit{Let $K\geq 1$. Let $\Omega$ be a domain in ${\mathbb{R}}^{m+n}$, \mbox{$m\geq n\geq 2$.}
Let $u:\, \Omega \rightarrow [-\infty ,+\infty )$ be  a Lebesgue measurable function. Suppose that the following conditions are satisfied:}
\begin{itemize}
\item[(a)] \textit{For each $y\in {\mathbb{R}}^n$ the function}
\[\Omega (y)\ni x\mapsto u(x,y)\in [-\infty ,+\infty )\]
\textit{is $K$-quasinearly subharmonic.}
 \item[(b)] \textit{For each $x\in {\mathbb{R}}^m$ the function}
\[\Omega (x)\ni y\mapsto u(x,y)\in [-\infty ,+\infty )\]
\textit{is $K$-quasinearly subharmonic.}
\item[(c)] \textit{There are strictly increasing surjections}  $\varphi :[0,+\infty )\rightarrow [0,+\infty )$ \textit{and} $\phi :[0,+\infty )\rightarrow [0,+\infty )$
 \textit{and}  \mbox{$s_0, \,s_1\in {\mathbb{N}}$,}  $s_0<s_1$, \textit{such  that}
\begin{itemize}
 \item[{(c1)}] \textit{$2K \phi ^{-1}(e^{s-s_0})\leq \phi ^{-1}(e^s)$ for all $s\geq s_1$,}
\item[{(c2)}] \textit{the following integral is convergent:}
\begin{equation*} \int\limits_{s_5}^{+\infty}\left[ \frac{\phi^{-1}(e^{s+2})}{\phi^{-1}(e^s)}
\frac{1}{\varphi (s-s_0)}\right]^{\frac{1}{m-1}}
 \left( \int\limits_{s_5+s_0+2}^{s+s_0+2}\left[ \frac{\phi^{-1}(e^{t+2})}{\phi^{-1}(e^t)} \right]^{\frac{1}{n-1}}\, dt
\right)^{\frac{n-1}{m-1}}\,ds<+\infty ,
\end{equation*}
\item[{(c3)}] $\varphi \circ \log^+ \phi (u^+)\in {\mathcal{L}}_{\textrm{loc}}^{1}(\Omega )$.
\end{itemize}
\end{itemize}
\textit{Then} $u$
\textit{is quasinearly subharmonic in $\Omega $}.

\vspace{2ex}

\noindent{\textbf{Corollary~4.3.}} {\textit{Let $K\geq 1$. Let $\Omega $ be a domain in ${\mathbb{R}}^{m+n}$, \mbox{$m\geq n\geq 2$.}
Let $u:\, \Omega \rightarrow [-\infty ,+\infty )$ be  a Lebesgue measurable function. Suppose that the following conditions are satisfied:}
\begin{itemize}
\item[(a)] \textit{For each $y\in {\mathbb{R}}^n$ the function}
\[\Omega (y)\ni x\mapsto u(x,y)\in [-\infty ,+\infty )\]
\textit{is $K$-quasinearly subharmonic.}
 \item[(b)] \textit{For each $x\in {\mathbb{R}}^m$ the function}
\[\Omega (x)\ni y\mapsto u(x,y)\in [-\infty ,+\infty )\]
\textit{is $K$-quasinearly subharmonic.}
\item[(c)] \textit{There are strictly increasing surjections}  $\varphi :[0,+\infty )\rightarrow [0,+\infty )$ \textit{and} $\phi :[0,+\infty )\rightarrow [0,+\infty )$
 \textit{and}  \mbox{$s_0, \,s_1\in {\mathbb{N}}$,}  $s_0<s_1$, \textit{such  that}
\begin{itemize}
 \item[{(c1)}] \textit{$2K \phi ^{-1}(e^{s-s_0})\leq \phi ^{-1}(e^s)$ for all $s\geq s_1$,}
\item[{(c2)}] $\phi^{-1}$ \textit{satisfies a $\Delta_2$-condition,}
\item[{(c3)}] \textit{the following integral is convergent:}
\begin{equation*} \int\limits_{s_1}^{+\infty}\frac{s^{\frac{n-1}{m-1}}}{{\varphi (s-s_0)}^{\frac{1}{m-1}}}\,ds<+\infty ,
\end{equation*}
\item[{(c4)}] $\varphi \circ \log^+ \phi ( u^+)\in {\mathcal{L}}_{\mathrm{loc}}^{1}(\Omega )$.
\end{itemize}
\end{itemize}
\textit{Then} $u$
\textit{is quasinearly subharmonic in $\Omega $}.

\vspace{2ex}

\noindent{\textit{Example~3.}}
 Let $u$ be separately subharmonic in $\Omega$. Let $\varphi :\, [0, +\infty )\rightarrow [0,+\infty )$ be a strictly increasing surjection such that 
\begin{equation*}\int\limits_{s_1}^{+\infty} \frac{s^{\frac{n-1}{m-1}}}
{\varphi (s-s_0)^{\frac{1}{m-1}}}ds<+\infty .
\end{equation*}
Choosing then various functions $\phi$, which, together with $\varphi$ and $u$, satisfy the conditions (c1), (c2) and (c4) of Corollary~4.3, one gets more concrete results. Possible choices are e.g.
\[\phi (t)=\frac{t^p}{(\log t)^q}, \, \, p>0, \, q\in {\mathbb{R}}.\]
 The case  $p=1$ and $q=0$ gives  the result of Armitage and Gardiner.

\vspace{2ex}

\noindent{\textit{Example~4.}}
 Let $u$ be separately subharmonic in $\Omega$ and $\varphi :\, [0, +\infty )\rightarrow [0,+\infty )$ be a strictly increasing surjection. Let $p>0$, $q\geq 0$.  Let $\phi :\, [0, +\infty )\rightarrow [0,+\infty )$ be a strictly increasing surjection for which there is $t_1>1$ such that 
\begin{equation*}
s=\phi (t)=e^{\left(\frac{\log t}{p}\right)^{\frac{1}{q+1}}}
=e^{\sqrt[q+1]{\frac{\log t}{p}}}      \,\, {\textrm{ for }}\,\, t\geq t_1,  
\end{equation*}
thus $t=\phi^{-1} (s)=e^{p(\log s)^{q+1}}$. 
One sees easily that the condition (c1) of Corollary~4.3 is satisfied, but the condition (c2) not. As a matter of fact, and as one easily sees,
\begin{equation*}
\frac{\phi^{-1}(e^{s+2})}{\phi^{-1}(e^{s})} \rightarrow +\infty \,\, {\textrm{ as }}\,\, s\rightarrow +\infty .
\end{equation*}  
Therefore, in this case one cannot use Corollary~4.3 to conclude that $u$ is subharmonic. However, using Corollary~4.2 we see that $u$ is subharmonic, provided 
that 
\begin{equation}\int\limits_{s_5}^{+\infty}\frac{e^{\frac{p}{m-1}[(s+2)^{q+1}-s^{q+1}]}}{\varphi (s-s_0)^{\frac{1}{m-1}}}\left( \int\limits_{s_5+s_0+2}^{s+s_0+2}e^{\frac{p}{n-1}[(t+2)^{q+1}-t^{q+1}]}dt\right)^{\frac{n-1}{m-1}} ds<+\infty ,
\end{equation}
and (this is just  (c3))
 \begin{equation*} 
\varphi \circ \left(\left[\frac{\log^+ u}{p}\right]^{\frac{1}{q+1}}\right)\in {\mathcal{L}}^1_{{\textrm{loc}}}(\Omega ).
\end{equation*}
The condition (4) is of course complicated, but it is easy to get simpler (but) stronger conditions, e.g. just estimating the inner integral. 

\vspace{2ex}
 
\noindent{\textit{Example~5.}}
 Let $u$ be separately subharmonic in $\Omega$ and $\varphi :\, [0, +\infty )\rightarrow [0,+\infty )$ be a strictly increasing surjection.
Let $p>0$ and  $\phi :\, [0, +\infty )\rightarrow [0,+\infty )$ be a strictly increasing surjection for which there is $t_1>1$ such that
\begin{equation*}
s=\phi (t)=( \log t)^p   \,\, {\textrm{ for }}\,\, t\geq t_1,  
\end{equation*}
and thus 
\begin{equation*}
t=\phi^{-1} (s)=e^{s^{\frac{1}{p}}}.  
\end{equation*}
Corollary~4.3 cannot now be applied, but from  Corollary~4.2 it follows that $u$ is subharmonic, provided that, in addition to the integrability condition (c3), 
\begin{equation*}
\varphi \circ \log^+ \circ ((\log^+ u^+)^p)\in {\mathcal{L}}^1_{{\textrm{loc}}}(\Omega ),
\end{equation*}
 also the condition (c2) holds. One possibility  to replace the, again rather complicated, condition (c2) by a simpler, but stronger one, is the following (we leave the details to the reader):
\begin{equation*}\int\limits_{s_5}^{+\infty}e^{\frac{2}{m-1}[e^{\frac{1}{p}(s+s_0+4)}-e^{\frac{1}{p}(s+s_0+2)}]} \frac{s^{\frac{n-1}{m-1}}}
{\varphi (s-s_0)^{\frac{1}{m-1}}}ds<+\infty .
\end{equation*}

\vspace{2ex}

\noindent{\textit{Remark~4.4.}
 As is seen above, the proof of Theorem~4.1 is essentially based on a previous simple result of separately quasinearly subharmonic functions, namely on  \cite{Rii07}, Proposition~3.1, p.~57, see also \cite{Rii09}, Proposition~3.2.1, p.~e2620. The situation is of course similar in the special case of separately subharmonic functions:  Armitage and Gardiner \cite{ArGa93}, proof of Theorem~1, pp.~257-259, base their result on the classical result of Avanissian \cite{Av61}, Théorème~9, p.~140. Equally well one might of course base the result on any  of the following later results: \cite{Le69}, Proposition~3, p.~24,   \cite{He71}, Theorem, p.~31,   Arsove \cite{Ar66}, Theorem~1, p.~622,  or  \cite{Rii89}, Theorem~1, p.~69. See also  Lelong \cite{Le45},  Théorème~1~bis, p.~315, and Cegrell and Sadullaev \cite{CeSa93}, Theorem~3.1, p.~82. Therefore, there are indeed good reasons to improve also these old basic results. In this connection, we  point out the following  recent improvement:  

\vspace{2ex}

\noindent{\textbf{Theorem~4.5.}}
 (\cite{Rii13}, Theorem~2, pp.~367-368) \textit{Let $K_1, K_2\geq 1$. Let $\Omega $ be a domain in ${\mathbb{R}}^{m+n}$, \mbox{$m,n\geq 2$.}
Let $u:\, \Omega \rightarrow
[-\infty ,+\infty )$ be such that}
 \begin{itemize}
\item[(a)] \textit{for each $y\in {\mathbb{R}}^n$ the function}
\[\Omega (y)\ni x\mapsto u(x,y)\in [-\infty ,+\infty )\]
\textit{is  $K_1$-quasinearly subharmonic, and, for almost every $y\in {\mathbb{R}}^n$, subharmonic,}
 \item[(b)] \textit{for each $x\in {\mathbb{R}}^m$ the function}
\[\Omega (x)\ni y\mapsto u(x,y)\in [-\infty ,+\infty )\]
\textit{is upper semicontinuous, and, for almost every $x\in {\mathbb{R}}^m$, \mbox{$K_2$-quasinearly subharmonic,}}
\item[(c)] \textit{for some $p>0$ there is a function}   $v\in {\mathcal{L}}_{\textrm{loc}}^p(\Omega )$ \textit{such that $u^+\leq v$}.
\end{itemize}
 \textit{Then for each $(a,b)\in \Omega $},
\begin{displaymath}\limsup_{(x,y)\rightarrow (a,b)}u(x,y)\leq K_1K_2\,u^+(a,b).\end{displaymath}

\vspace{2ex}

Observe that the proof of the above (quasinearly subharmonicity) result, and thus also the proof of the following special case result, is simpler than the proofs of the older subharmonicity results. (See also the previous versions \cite{Rii07}, Corollary~3.2 and Corollary~3.3, p.~61, and  \cite{Rii09}, Corollary~3.2.4 and Corollary~3.2.5, p.~e2621).

\vspace{2ex}

\noindent{\textbf{Corollary~4.6.}}
 (\cite{Rii13}, Corollary~2, p.~369) \textit{ Let $\Omega $ be a domain in ${\mathbb{R}}^{m+n}$, \mbox{$m,n\geq 2$.}
Let $u:\, \Omega \rightarrow
[-\infty ,+\infty )$ be such that}
 \begin{itemize}
\item[(a)] \textit{for each $y\in {\mathbb{R}}^n$ the function}
\[\Omega (y)\ni x\mapsto u(x,y)\in [-\infty ,+\infty )\]
\textit{is nearly subharmonic, and, for almost every $y\in {\mathbb{R}}^n$, subharmonic,}
\item[(b)] \textit{for each} $x\in {\mathbb{R}}^m$ \textit{the function}
\[\Omega (x)\ni y\mapsto u(x,y)\in [-\infty ,+\infty )\]
\textit{is upper semicontinuous, and, for almost every $x\in {\mathbb{R}}^m$, (nearly) subharmonic,} 
\item[(c)] \textit{for some $p>0$ there is a function}   $v\in {\mathcal{L}}_{\mathrm{loc}}^p(\Omega )$ \textit{such that $u^+\leq v$}.
\end{itemize}
 \textit{Then} $u$
\textit{is  upper semicontinuous and thus subharmonic in $\Omega$}.

\vspace{2ex}

\noindent{\textit{Proof}}. It is easy to see that for each $M\geq 0$, the function  $u_M:=\max\{u,-M\}+M$ satisfies the assumptions of Theorem~4.5. Thus $u_M$ is upper semicontinuous. Since by \cite{Rii07}, Corollary~3.1, p.~59 (see also \cite{Rii09}, Corollary~3.2.3, pp.~e2620-e2621, and \cite{Rii13}, Corollary~1, p.~367),  $u_M$ is anyway nearly subharmonic, it is in fact subharmonic. Using then a standard result, see e.g. \cite{He71}, a), p.~8, one sees that $u$ is subharmonic and thus also upper semicontinuous. \qed

\vspace{2ex}

\noindent{\textit{Remark~4.7.}} Observe that Corollary~4.6 is partially related to the result \cite{He71}, Proposition~2, pp.~34--35. Though our assumptions are partly slightly stronger, our proof (see \cite{Rii13}, pp.~367-370) is, on the other hand,  different and shorter.

\end{document}